\documentclass[twoside,a4paper]{amsart}

\usepackage{xcolor}
\definecolor{webgreen}{rgb}{0,.5,0}
\definecolor{webbrown}{rgb}{.6,0,0}
\definecolor{RoyalBlue}{cmyk}{1, 0.50, 0, 0}
\usepackage[colorlinks=true, breaklinks=true, urlcolor=webbrown, linkcolor=RoyalBlue, citecolor=webgreen,backref=page]{hyperref}
\usepackage{epsfig, graphicx, subfigure}
\usepackage{verbatim, setspace,pifont}
\usepackage{amsmath, amssymb}

\newcommand{\N}{{\mathbb N}}

\newcommand{\R}{{\mathbb R}}
\newcommand{\C}{{\mathbb C}}
\newcommand{\D}{{\mathbb D}}
\newcommand{\T}{{\mathbb T}}

\newcommand{\E}{{\mathbb E}}

\newcommand{\re}{{\mathsf{Re}}}
\newcommand{\im}{{\mathsf{Im}}}

\newcommand{\dist}{{\mathsf{dist}}}

\newcommand{\dd}{{\mathrm d}}
\newcommand{\ic}{{\mathrm i}}

\newtheorem{lemma}{Lemma}[section]
\newtheorem{theorem}[lemma]{Theorem}

\newtheorem{proposition}[lemma]{Proposition}
\newtheorem{corollary}[lemma]{Corollary}

\begin{document}

\title[Zeros of Real Random Polynomials Spanned by OPUC]{Zeros of Real Random Polynomials \\ Spanned by OPUC}

\author{Maxim L. Yattselev}
\address{Department of Mathematical Sciences, Indiana University-Purdue University Indianapolis, 402~North Blackford Street, Indianapolis, IN 46202, USA}
\email{maxyatts@iupui.edu}

\author{Aaron Yeager}
\address{Department of Mathematics, Oklahoma State University, Stillwater, OK 74078 USA}
\email{aaron.yeager@okstate.edu}

\thanks{The research of M. Yattselev was supported by a grant from the Simons Foundation, CGM-354538.}

\begin{abstract}
Let \( \{\varphi_i\}_{i=0}^\infty \) be a sequence of orthonormal polynomials on the unit circle with respect to a probability measure \( \mu \). We study zero distribution of random linear combinations of the form
\[
P_n(z)=\sum_{i=0}^{n-1}\eta_i\varphi_i(z),
\]
where \( \eta_0,\dots,\eta_{n-1} \) are i.i.d. standard Gaussian variables.  We use the Christoffel-Darboux formula to simplify the density functions provided by Vanderbei for the expected number real and complex of zeros  of \( P_n \).  From these expressions, under the assumption that \( \mu \) is in the Nevai class, we deduce the limiting value of these density functions away from the unit circle.  Under the mere assumption that \( \mu \) is doubling on subarcs of \( \T \) centered at \( 1 \) and \( -1 \), we show that the expected number of real zeros of \( P_n \) is at most
\[
(2/\pi) \log n +O(1),
\]
and that the asymptotic equality holds when the corresponding recurrence coefficients decay no slower than \( n^{-(3+\epsilon)/2} \), \( \epsilon>0 \).  We conclude with providing results that estimate the expected number of complex zeros of \( P_n \) in shrinking neighborhoods of compact subsets of \( \T \).
\end{abstract}

\keywords{}

\subjclass[2010]{Primary: 30C15; Secondary: 30B20, 60B10}

\maketitle

\section{Introduction and Main Results}
\label{s:intro}

In \cite{Kac43}, Kac considered random polynomials
\begin{equation}
\label{Kac}
P_n(z) = \eta_0 + \eta_1 z + \cdots +\eta_{n-1}z^{n-1},
\end{equation}
where \( \eta_i \) are i.i.d. standard real Gaussian variables. He has shown that \( \E(N_n(\Omega)) \), the expected number of zeros of \( P_n \) on a measurable set \( \Omega\subset\R \), is equal to
\begin{equation}
\label{real-intensity-Kac}
\E(N_n(\Omega)) = \frac1\pi\int_\Omega\frac{\sqrt{1-h_n^2(x)}}{|1-x^2|}\dd x, \quad h_n(x) = \frac{nx^{n-1}(1-x^2)}{1-x^{2n}},
\end{equation}
from which he proceeded with an estimate
\begin{equation}
\label{real-expected-Kac}
\E(N_n(\R)) = \frac{2+o(1)}\pi\log n \quad \text{as} \quad n\to\infty.
\end{equation}
In fact, it was eventually shown by Wilkins \cite{Wilk88} that \( \E(N_n(\R)) \) has an asymptotic expansion of the form
\[
\E(N_n(\R)) \sim \frac2\pi\log n + \sum_{k=0}^\infty A_kn^{-k}
\]
for some constants \( A_k \). In another connection, Shepp and Vanderbei \cite{ShVan95} showed that
\begin{equation}
\label{expected-SV}
\E(N_n(\Omega)) = \int_{\Omega\cap\R} \rho_n^{(1,0)}(x)\mathrm dx + \int_\Omega\rho_n^{(0,1)}(z)\mathrm dx\mathrm dy,
\end{equation}
where \( \Omega \subset \C \) is measurable, \( \rho_n^{(1,0)}(x) \) is the integrand in \eqref{real-intensity-Kac}, and \( \rho_n^{(0,1)}(z) \) is the explicitly known intensity function for the complex zeros, see \eqref{real-intensity-Vanderbei} and \eqref{complex-intensity-Vanderbei} further below. At the same time Edelman and Kostlan \cite{EdKos95} considered random functions of the form
\begin{equation}
\label{Vanderbei}
P_n(z) = \eta_0f_0(z) + \eta_1 f_1(z) + \cdots +\eta_{n-1}f_{n-1}(z),
\end{equation}
where \( \eta_i \) are certain real random variables and \( f_i(z) \) are arbitrary functions on the complex plane that are real on the real line. Using beautiful and simple geometrical argument they have shown\footnote{In fact, Edelman and Kostlan derive an expression for the real intensity function for any random vector \( (\eta_0,\ldots,\eta_{n-1}) \) in terms of its joint probability density function and of \( v(x) \).} that if \( \eta_0,\ldots,\eta_{n-1} \) are elements of a multivariate real normal distribution with mean zero and covariance matrix \( C \) and the functions \( f_i \) are differentiable on the real line, then
\begin{equation}
\label{Edelman-Kostlan}
\rho_n^{(1,0)}(x) = \left.\frac1\pi\frac{\partial^2}{\partial s\partial t}\log\left( v(s)^\mathsf{T}Cv(t)\right)\right|_{t=s=x}, \end{equation}
where \( v(x) = \big(f_0(x),\ldots,f_{n-1}(x)\big)^\mathsf{T} \). If random variables \( \eta_i \) in \eqref{Vanderbei} are again i.i.d. standard real Gaussians, \eqref{Edelman-Kostlan} specializes to
\begin{equation}
\label{real-intensity-Vanderbei}
\rho_n^{(1,0)}(x) = \frac1\pi\frac{\sqrt{K_n(x,x)K_n^{(1,1)}(x,x)-K_n^{(1,0)}(x,x)^2}}{K_n(x,x)}
\end{equation}
(\eqref{real-intensity-Vanderbei} was also independently rederived in \cite[Proposition~1.1]{LubPritXie17} and \cite[Theorem~1.2]{uVan}), where
\begin{equation}
\label{kernels}
\left\{
\begin{array}{lll}
K_n(z,w) &:=& \sum_{i=0}^{n-1}f_i(z)\overline{f_i(w)}, \medskip \\
K_n^{(1,0)}(z,w) & := & \sum_{i=0}^{n-1}f_i^\prime(z)\overline{f_i(w)}, \medskip \\
K_n^{(1,1)}(z,w) &:=& \sum_{i=0}^{n-1}f_i^\prime(z)\overline{f_i^\prime(w)}.
\end{array}
\right.
\end{equation}
Moreover, if the functions \( f_i \) in \eqref{Vanderbei} are in addition entire, then it was shown by Vanderbei \cite[Theorem~1.1]{uVan} that
\begin{equation}
\label{complex-intensity-Vanderbei}
\begin{array}{lll}
\rho_n^{(0,1)}(z) &=& \displaystyle \frac1\pi\frac{K_n^{(1,1)}(z,z)}{\big(K_n(z,z)^2-|K_n(z,\overline z)|^2\big)^{1/2}} \medskip \\
&-& \displaystyle \frac1\pi\frac{K_n(z,z)\big(|K_n^{(1,0)}(z,z)|^2+|K_n^{(1,0)}(z,\overline z)|^2\big)}{\big(K_n(z,z)^2-|K_n(z,\overline z)|^2\big)^{3/2}} \medskip \\
&+& \displaystyle \frac2\pi\frac{\re\big(K_n(z,\overline z)K_n^{(1,0)}(z,z)K_n^{(1,0)}(\overline z,z)\big)}{\big(K_n(z,z)^2-|K_n(z,\overline z)|^2\big)^{3/2}}.
\end{array}
\end{equation}

The result of Kac was generalized in many directions, see \cite{EdKos95} and the introduction of \cite{uVan} for the references. In this note we concentrate on the case where the functions \( f_i=\varphi_i \) are real orthonormal polynomials on the unit circle complementing the case where \( f_i=p_i \) are orthonormal polynomials on the real line \cite{LubPritXie17}. That is, for some probability Borel measure \( \mu \) on \( \T \) that is symmetric with respect to conjugation, it holds that
\begin{equation}
\label{ortho}
\int_\T \varphi_i(z)\overline{\varphi_j(z)}\dd\mu(z) = \delta_{ij},
\end{equation}
where \( \delta_{ij} \) is the usual Kronecker symbol. Notice that all the coefficients of the polynomials \( \varphi_i \) are real due to the conjugate-symmetry of the measure \( \mu \). Write \( \varphi_i(z)=\kappa_i\Phi_i(z) \), where \( \Phi_i \) is monic. Polynomials \( \Phi_i \) satisfy Szeg\H{o} recurrence relations
\begin{equation}
\label{recurrence}
\left\{
\begin{array}{l}
\displaystyle \Phi_{i+1}(z) = z\Phi_i(z) - \alpha_i\Phi_i^*(z), \medskip \\
\displaystyle \Phi_{i+1}^*(z) = \Phi_i^*(z) - \alpha_iz\Phi_i(z),
\end{array}
\right.
\end{equation}
where \( \Phi_i^*(z):=z^i\Phi_i(1/z) \) and \( \{\alpha_i\}\subset(-1,1) \) are the recurrence coefficients, see \cite[Theorem~1.5.2]{Simon1}. Recall also that given \( \{\alpha_i\}\subset(-1,1) \), there exists a unique conjugate-symmetric probability measure \( \mu \) whose monic orthogonal polynomials satisfy \eqref{recurrence}, see \cite[Theorem~1.7.11]{Simon1}. Moreover, in this case it holds that
\begin{equation}
\label{leading-coeff}
\kappa_n = \prod_{i=0}^{n-1}\big(1-\alpha_i^2\big)^{-1/2},
\end{equation}
see \cite[Equation~(1.5.22)]{Simon1}. Thus, we might refer to orthonormal polynomials \( \varphi_i \) as defined by either the measure \( \mu \) or the recurrence coefficients \( \{\alpha_i\} \).

\begin{theorem}
\label{thm:real-intensity}
Let \( \{\varphi_i \} \) be a sequence of polynomials satisfying \eqref{ortho}. Further, let \( P_n \) be a real random polynomial \eqref{Vanderbei} with \( f_i=\varphi_i \) and \( \eta_i \) being i.i.d. standard real Gaussian variables. Then the intensity function \( \rho_n^{(1,0)} \) from \eqref{real-intensity-Vanderbei} can be written as
\begin{equation}
\label{OPUC-real-intensity}
 \rho_n^{(1,0)}(x) = \frac{1}{\pi}\frac{\sqrt{1-h_n^2(x)}}{|1-x^2|}, \quad h_n(x) = \frac{(1-x^2)b_n^\prime(x)}{1-b_n^2(x)}, \quad b_n(z) = \frac{\varphi_n(z)}{\varphi_n^*(z)}.
\end{equation}
\end{theorem}

Clearly, if the recurrence coefficients are all zero, \( \varphi_n(z) = z^n \) and respectively \( b_n(z) = z^n \). That is, we recover the intensity function from the Kac formula \eqref{real-intensity-Kac}.

\begin{proposition}
\label{prop:Caratheodory}
Under the conditions of Theorem~\ref{thm:real-intensity}, it holds that
\[
h_n(z) = \frac{1-z^2}2\frac{F_n^\prime(z)}{F_n(z)},
\]
where \( F_n \) is a trivial Caratheodory function given by
\[
F_n(z) = -\frac{\Phi_n(z;-1)}{\Phi_n(z;1)} = \int\frac{\zeta+z}{\zeta-z}\dd\sigma_n(\zeta),
\]
\( \sigma_n \) is a probability measure defined by
\[
\sigma_n = \sum_{k=1}^n \frac{|\varphi_{n-1}(\zeta_{k,n})|^2}{K_n(\zeta_{k,n},\zeta_{k,n})}\delta_{\zeta_{k,n}},
\]
\( \delta_\zeta \) is the unit point mass at \( \zeta \), \( \zeta_{k,n}\in \T \) are the zeros of \( \Phi_n(z;1) \), and
\[
\Phi_n(z;\beta) := z\Phi_{n-1}(z) - \beta\Phi_{n-1}^*(z), \quad |\beta|=1,
\]
is the \( n \)-th paraorthogonal polynomial\footnote{It is known \cite[Theorem~2.2.12]{Simon1} that all the zeros of \( \Phi_n(z;\beta) \) are simple and lie on the unit circle.} associated with parameter \( \beta \).
\end{proposition}

Since the Blaschke products \( b_n \) necessarily satisfy \( |b_n(z)|\leq 1 \) in \( \overline \D \), they form a normal family there. Moreover, as \( b_n(1/z) =1/b_n(z) \), we see that
\begin{equation}
\label{symmetry}
\rho_n^{(1,0)}(1/x) = x^2\rho_n^{(1,0)}(x).
\end{equation}
The following corollary is immediate.

\begin{corollary}
In the setting of Theorem~\ref{thm:real-intensity}, let \( \mathcal N\subset \N \) be such that \( b_n(z)\to b(z)\not\equiv1 \) as \( \mathcal N\ni n\to\infty \) for some analytic function \( b(z) \) in \( \D \). Then
\[
\rho_n^{(1,0)}(x) \to \frac1\pi\frac{\sqrt{1-h^2(x)}}{1-x^2}, \quad h(x) = b^\prime(x)\frac{1-x^2}{1-b^2(x)},
\]
locally uniformly on \( (-1,1) \) as \( \mathcal N\ni n\to\infty \). In particular, if \( \alpha_i \to 0 \) as  \( i\to\infty \), then
\[
\rho_n^{(1,0)}(x) \to \frac1\pi\frac1{|1-x^2|} \quad \text{as} \quad n\to\infty
\]
uniformly on closed subsets of \( (-\infty,-1)\cup(-1,1)\cup(1,\infty) \) as \( n\to\infty \).
\end{corollary}

The second claim of the corollary is a straightforward consequence of the fact that \( b_n(z)\to 0 \) locally uniformly in \( \D \) as \( n\to\infty \) if and only if \( \alpha_i \to 0 \) as  \( i\to\infty \), see \cite[Theorem~1.7.4]{Simon1}.

Recall that a measure \( \mu \) is called doubling on a subarc \( T \subseteq \T \) if there exists a constant \( L> 0 \) such that
\[
\mu(2I) \leq L \mu(I), \quad 2I\subseteq T,
\]
for any subarc \( I \), where \( 2I \) is a subarc of \( T \) with the same center as \( I \) and twice the arclength.

\begin{theorem}
\label{thm:real-expected}
In the setting of Theorem~\ref{thm:real-intensity}, assume that there exist two subarcs of \( \T \), centered at \( 1 \) and \( -1 \), on which \( \mu \) is doubling. Then it holds that
\begin{equation}
\label{upper-bound}
\E(N_n(\R)) \leq \frac2\pi\log n + \mathcal O(1).
\end{equation}
Moreover, if the quantities \( |k^p\alpha_k| \) are uniformly bounded above for some \( p>3/2 \), then
\begin{equation}
\label{exact-asymp}
\E(N_n(\R)) = \frac{2+o(1)}\pi\log n.
\end{equation}
\end{theorem}

The assumption on \( \alpha_k \)'s in \eqref{exact-asymp} implies that they are absolutely summable. Hence, it follows from Baxter's theorem, see \cite[Theorem~5.2.1]{Simon1}, that \( \mu \) is absolutely continuous with respect to the arclenth distribution on \( \T \) and the Radon-Nikodym derivative is continuous and non-vanishing there. In particular, \( \mu \) is doubling on \( \T \).

\begin{proposition}
\label{prop:real-expected}
In the setting of Theorem~\ref{thm:real-intensity}, assume that
\[
\mu = t\nu + (1-t)\delta_1, \quad t\in(0,1),
\]
where \( \nu \) is a conjugate-symmetric probability measure on the unit circle such that  \( |k^p\alpha_k(\nu)| \) are uniformly bounded above for some \( p>3/2 \). Then \eqref{exact-asymp} holds while
\begin{equation}
\label{new-rec}
\alpha_{n-1} = \alpha_{n-1}(\nu) + \varphi_{n-1}(1;\nu)\varphi_n(1;\nu)\frac{\sqrt{1-|\alpha_{n-1}(\nu)|^2}}{t(1-t)^{-1}+K_n(1,1;\nu)},
\end{equation}
where the quantities \( \alpha_n(\nu), \varphi_n(z;\nu), K_n(z,w;\nu) \) are defined as before only with respect to the measure \( \nu \). 
\end{proposition}

Formula \eqref{new-rec} was derived in \cite{Wong09} and shows that \(\alpha_n\sim 1/n \) as \( n\to\infty \), that is, the recurrence coefficients do not obey the conditions of Theorem~\ref{thm:real-expected}. Indeed, the coefficients \( \alpha_k(\nu) \) are absolutely summable. Thus, there exists a constant \( c>1 \) such that \( c^{-1} \leq |\varphi_k(1;\nu)| \leq c \) for all \( k \), see \cite[Equation (1.5.16)]{Simon1} and \eqref{leading-coeff}. Hence, \( n/c^2 \leq K_n(1,1;\nu) \leq nc^2 \), which yields the claim.

\begin{theorem}
\label{thm:complex-intensity}
In the setting of Theorem~\ref{thm:real-intensity}, assume that \( \alpha_k \to 0 \) as  \( k\to\infty \). Then we have
\[
\rho_n^{(0,1)}(z) \to \frac{1}{\pi(1-|z|^2)^2}\sqrt{1-\left|\frac{1-|z|^2}{1-z^2}\right|^2   }
\]
locally uniformly in \( \overline\C\setminus(\T\cup\R) \) as \( n\to\infty \).
\end{theorem}

It follows from Theorem~\ref{thm:complex-intensity} that the zeros of \( P_n \) almost surely accumulate on the unit circle. In fact, the following discrepancy results hold.

\begin{theorem}
\label{thm:complex-expected}
In the setting of Theorem~\ref{thm:real-intensity}, assume that the measure of orthogonality \( \mu \) in \eqref{ortho} is regular in the sense of Ullman-Stahl-Totik, that is,
\begin{equation}
\label{Reg}
\varepsilon_n := \frac1n\log|\kappa_n| \to 0 \quad \text{as} \quad n\to\infty,
\end{equation}
where \( \kappa_n \) is the leading coefficient of \( \varphi_n \). Then for any subarc \( S \subset \T \), it holds that
\begin{equation}
\label{discrepancy}
\E\left[\left|\frac1nN_n\left(\Omega\left(S,\delta\right)\right) -  \frac{|S|}{2\pi}\right|\right] \leq \frac C\delta \sqrt{\frac{\log n}{n} +\sqrt{\varepsilon_n}  } \quad \text{as} \quad n\to\infty,
\end{equation}
where \( C\) is independent of \( n \) and \( \Omega(S,\delta) := \big\{ rz: z\in S, \; r\in(1-\delta,1+\delta) \big\} \).
\end{theorem}

The claim of Theorem~\ref{thm:complex-expected} is contained in \cite[Corollary~3.2]{PritYgr15} under the additional conditions that the measure $\mu$ for the basis \( \{\varphi_i\} \) is absolutely continuous with Radon-Nikodym derivative bounded below by a positive constant on the whole unit circle. We show that the proof of \cite[Corollary~3.2]{PritYgr15} in fact remains valid under restriction \eqref{Reg}.

It follows from \eqref{discrepancy} that the zeros of \( P_n \) that are expected to approach an arc \( S \) are contained in an annular neighborhood of width \( \big(\frac{\log n}n+\sqrt{\epsilon_n}\big)^{1/2-\epsilon} \) for any \( \epsilon>0 \). More stringent assumptions on the measure \( \mu \) allow us to decrease the width this neighborhood.

\begin{theorem}
\label{thm:complex-expected2}
In the setting of Theorem~\ref{thm:real-intensity}, assume that \( \alpha_i \to 0 \) as  \( i\to\infty \). Let \( S \) be a compact subset of \( \T\setminus\{\pm1\} \). Assume, in addition, that \( \mu \) is absolutely continuous with respect to the arclength measure on an open set containing \( S \) and its Radon-Nikodym derivative is positive and continuous at each point of \( S \). Given \( -\infty<\tau_1<\tau_2<\infty \), it follows that
\begin{equation}
\label{discrepancy2}
\frac1n \E\big[N_n\big(\Omega(S,\tau_1,\tau_2)\big)\big] \to \frac{|S|}{2\pi}\left(\frac{H^\prime(\tau_2)}{H(\tau_2)}-\frac{H^\prime(\tau_1)}{H(\tau_1)}\right) \quad \text{as} \quad n\to\infty,
\end{equation}
where \( \Omega(S,\tau_1,\tau_2) := \big\{ rz: z\in S, \; r\in(1+\frac{\tau_1}{2n},1+\frac{\tau_2}{2n}) \big\} \) and \( \displaystyle H(\tau) := \frac{e^\tau-1}{\tau} \).
\end{theorem}

It can be readily verified that \( H^\prime/H \) is increasing on the real line with
\[
\lim_{\tau\to-\infty}\frac{H^\prime(\tau)}{H(\tau)}=0 \quad \text{and} \quad \frac{H^\prime(\tau)}{H(\tau)} =  1 - \frac{H^\prime(-\tau)}{H(-\tau)}.
\]
Hence, the zeros of \( P_n \) approaching \( S \) are expected to be contained in an annular band around \( S \) of width \( n^{-1+\epsilon} \) for any \( \epsilon>0 \).

\section{Proof of Theorem~\ref{thm:real-intensity} and Proposition~\ref{prop:Caratheodory}}

\subsection{Proof of Theorem~\ref{thm:real-intensity}}

According to the Christoffel-Darboux formula \cite[Theorem~2.2.7]{Simon1} and since polynomials \( \varphi_n \) have real coefficients, it holds that
\begin{equation}
\label{Kn}
K_n(z,w) = \frac{\varphi_n^*(z)\varphi_n^*(\overline w) -\varphi_n(z)\varphi_n(\overline w)}{1-z\overline w}.
\end{equation}
Hence,
\begin{equation}
\label{Kn1}
K_n^{(1,0)}(z,w) = \frac{(\varphi_n^*)^\prime(z)\varphi_n^*(\overline w)-\varphi_n^\prime(z)\varphi_n(\overline w)}{1-z\overline w} + \overline w\frac{K_n(z,w)}{1-z\overline w}
\end{equation}
and
\begin{multline}
\label{Kn11}
K_n^{(1,1)}(z,w) = \frac{(\varphi_n^*)^\prime(z)(\varphi_n^*)^\prime(\overline w)-\varphi_n^\prime(z)\varphi_n^\prime(\overline w)}{1-z\overline w} + z\frac{(\varphi_n^*)^\prime(z)\varphi_n^*(\overline w)-\varphi_n^\prime(z)\varphi_n(\overline w)}{(1-z\overline w)^2} \\ + \overline w \frac{\varphi_n^*(z)(\varphi_n^*)^\prime(\overline w) -\varphi_n(z)\varphi_n^\prime(\overline w)}{(1-z\overline w)^2} + \frac{(1+z\overline w)K_n(z,w)}{(1-z\overline w)^2}.
\end{multline}
Thus,
\[
K_n(x,x)K_n^{(1,1)}(x,x)-K_n^{(1,0)}(x,x)^2 = \frac{K_n^2(x,x)}{(1-x^2)^2} - \left(\frac{\varphi_n^*(x)\varphi_n^\prime(x)-\varphi_n(x)(\varphi_n^*)^\prime(x)}{1-x^2}\right)^2.
\]
Therefore, the claim of the theorem now follows from \eqref{real-intensity-Vanderbei} since
\[
\rho_n^{(1,0)}(x)^2 = \frac1{\pi^2}\left[\frac1{(1-x^2)^2} - \left(\frac{b_n^\prime(x)}{1-b_n^2(x)}\right)^2\right] = \frac1{\pi^2}\frac{1-h_n^2(x)}{(1-x^2)^2}.
\]

\subsection{Proof of Proposition~\ref{prop:Caratheodory}}

It can be readily verified using recurrence relations \eqref{recurrence} that
\begin{equation}
\label{hn-rep}
h_n(z) = (1-z^2)\frac{(zb_{n-1}(z))^\prime}{1-(zb_{n-1}(z))^2}.
\end{equation}
Therefore, the first claim of the proposition holds with
\[
F_n(z) := \frac{1+zb_{n-1}(z)}{1-zb_{n-1}(z)}.
\]
Since \( zb_{n-1}(z) \) is a finite Blaschke product, \( F_n(z) \) is a trivial Caratheodory function \cite[Section~1.3]{Simon1}. Clearly,
\[
F_n(z) = -\frac{z\Phi_{n-1}(z)+\Phi_{n-1}^*(z)}{z\Phi_{n-1}(z)-\Phi_{n-1}^*(z)} = - \frac{\Phi_n(z;-1)}{\Phi_n(z;1)}
\]
by the very definition of \( b_{n-1}(z) \) and \( \Phi_n(z;\beta) \). Thus, it only remains to prove the integral representation. On the one hand, it holds that
\[
\int_\T \frac{\zeta+z}{\zeta-z}\dd\sigma_n = -\sum_{k=1}^n \frac{|\varphi_{n-1}(\zeta_{k,n})|^2}{K_n(\zeta_{k,n},\zeta_{k,n})}\frac{z+\zeta_{k,n}}{z-\zeta_{k,n}} = -\left(\sum_{k=1}^n\frac{|\varphi_{n-1}(\zeta_{k,n})|^2}{K_n(\zeta_{k,n},\zeta_{k,n})}\right)\frac{z^n+\cdots}{z^n+\cdots},
\]
which is a rational function of type \( (n,n) \) with a simple pole at each \( \zeta_{k,n} \) with residue
\[
-2\zeta_{k,n}\frac{|\varphi_{n-1}(\zeta_{k,n})|^2}{K_n(\zeta_{k,n},\zeta_{k,n})}.
\]
On the other hand, by the very definition, \( F_n(z) \) is also a rational function of type \( (n,n) \) with a simple pole at each \( \zeta_{k,n} \). Next observe that
\begin{equation}
\label{para}
\Phi_n(z;1) = \frac{z-\zeta_{k,n}}{\kappa_{n-1}\overline{\varphi_{n-1}(\zeta_{n,k})}}K_n(z,\zeta_{k,n}),
\end{equation}
where \( \kappa_{n-1} \) is the leading coefficient of \( \varphi_{n-1}(z) \). Indeed, the normalization follows from the fact that the left-hand side is a monic polynomial. Furthermore, Christoffel-Darboux formula \eqref{Kn} gives
\begin{multline}\label{para2}
\kappa_n^{-2}\overline\zeta(\zeta-z)K_n(z,\zeta) = \Phi_n^*(z)\overline{\Phi_n^*(\zeta)} - \Phi_n(z)\overline{\Phi_n(\zeta)} \\
 = \Phi_{n-1}^*(z)\big(\overline{\Phi_n^*(\zeta) + \alpha_{n-1}\Phi_n(\zeta)}\big) - z\Phi_{n-1}(z)\big(\overline{\alpha_{n-1}\Phi_n^*(\zeta) + \Phi_n(\zeta)}\big)
\end{multline}
for any \( \zeta\in\T \), where we used \eqref{recurrence} for the second equality. Now, since
\begin{equation}
\label{zero}
\zeta_{k,n}\Phi_{n-1}(\zeta_{k,n}) = \Phi_{n-1}^*(\zeta_{k,n}),
\end{equation}
we get from \eqref{recurrence} that
\begin{equation*}
\Phi_n^*(\zeta_{k,n}) + \alpha_{n-1}\Phi_n(\zeta_{k,n}) = \Phi_{n-1}^*(\zeta_{k,n})(1 - \alpha_{n-1}^2) = \alpha_{n-1}\Phi_n^*(\zeta_{k,n}) + \Phi_n(\zeta_{k,n}).
\end{equation*}
Combining the above equality with \eqref{para2} finishes the justification of \eqref{para}. Now, \eqref{para} and \eqref{zero} yield that the residue of \( F_n \) at \( \zeta_{k,n} \) is equal to
\[
-2\zeta_{k,n}\Phi_{n-1}(\zeta_{k,n})\frac{\kappa_{n-1}\overline{\varphi_{n-1}(\zeta_{n,k})}}{K_n(\zeta_{k,n},\zeta_{k,n})} = -2\zeta_{k,n}\frac{|\varphi_{n-1}(\zeta_{n,k})|^2}{K_n(\zeta_{k,n},\zeta_{k,n})},
\]
which proves the integral representation. In particular, since \( F_n(z) = -1 + O(1/z) \) as \( z\to\infty \), it must hold that
\[
\sum_{k=1}^n\frac{|\varphi_{n-1}(\zeta_{k,n})|^2}{K_n(\zeta_{k,n},\zeta_{k,n})} = 1,
\]
that is, \( \sigma_n \) is a probability measure.

\section{Proof of Theorem~\ref{thm:real-expected} and Proposition~\ref{prop:real-expected}}

In what follows, to avoid complicated schemes of labeling constants, we shall write
\[
f_n(z) \lesssim g_n(z), \quad z\in K,~n\in\N \quad \Leftrightarrow \quad f_n(z) \leq Cg_n(z), \quad z\in K,~n\in\N,
\]
where the constant \( C \) depends possibly on \( K \) but not on \( z \). Furthermore, we write
\[
f_n(z) \sim g_n(z)  \quad \Leftrightarrow \quad f_n(z) \lesssim g_n(z) \lesssim f_n(z).
\]

\subsection{An Auxiliary Estimate}

Recall that the \( n \)-th Christoffel function of \( \mu \) is given by
\begin{equation}
\label{Christoffel}
\lambda_n(z;\mu) := \inf_{\deg(p)\leq n-1}|p(z)|^{-2}\int|p|^2\dd \mu = K_n^{-1}(z,z),
\end{equation}
where the last equality is extremely well known, see for example \cite[Equation~(1.2.39)]{Simon1}. In this subsection we prove the following claim: \emph{if the measure \( \mu \) is doubling on a subarc \( T\subset\T \), then it holds that}
\begin{equation}
\label{Christoffel-bounds}
\lambda_n\big(ze^{\ic a/n};\mu\big) \sim \mu_n(z) := \int_{T(z,\frac1n)}\dd\mu, \quad z\in T^\prime, \quad |a|\leq2,
\end{equation}
uniformly with respect to \(z,a,n \), where \( T^\prime\subset T \) is a subarc with endpoints different from those of \( T \) and \( T(z,\delta) \) stands for the subarc of \( \T \) centered at \( z \) of arclength \( 2\delta \). When \( \mu \) is doubling on the whole circle \( \T \) and \( a=0 \), this claim is simply \cite[Theorem~4.3]{MastTot00}. The proof of the localized version \eqref{Christoffel-bounds} is quite similar to the one of \cite[Theorem~4.3]{MastTot00}. However, to improve readability, we adduce the full proof of this fact below.

Given an integer \( m \) that we shall fix later, put
\begin{eqnarray*}
S_n(z,\eta) & := & \gamma_n \left(\sum_{k=0}^{n-1}\bigg(\frac z\eta\bigg)^k\right)^m \left(\sum_{k=0}^{n-1}\bigg(\frac\eta z\bigg)^k\right)^m \\
& = & \gamma_n\left(\frac1{z\eta}\right)^{m(n-1)}\left(\frac{z^n-\eta^n}{z-\eta}\right)^{2m} = \gamma_n\left(\frac{\sin(\frac{a-b}2)}{\sin(\frac{a-b}{2n})}\right)^{2m}
\end{eqnarray*}
where \( z=e^{\ic a/n} \), \( \eta=e^{\ic b/n} \), and the normalizing constant \( \gamma_n \) is chosen so that
\(
\int_\T S_n(z,\eta)|\dd \eta| =1.
\)
It is known that \( \gamma_n\sim n^{-2m+1} \). The last representation of \( S_n(z,\eta) \) shows that
\begin{equation}
\label{Sn-lower}
\re\left(S_n\big(e^{\ic a/n},e^{\ic b/n}\big)\right) \gtrsim n.
\end{equation}
locally uniformly for \( |a-b|<2\pi \). Similarly, we can easily see from the second representation that
\begin{equation}
\label{Sn-upper}
|S_n(z,\eta)| \lesssim \left\{
\begin{array}{ll}
n, & |z-\eta|\leq \frac1n, \medskip \\
n^{-2m+1}|z-\eta|^{-2m}, & |z-\eta|\geq \frac1n,
\end{array}
\right.
\end{equation}
for \( |n(|z|-1)|,|n(|\eta|-1)|\leq A \), where the constant is uniform in \( A>0 \).

We start with an upper bound. Let \( z\in T^\prime \) and \( |a|\leq 2 \). Since \( S_n\big(ze^{\ic a},ze^{\ic a}\big) =\gamma_nn^{2m} \sim n \), it is immediate from \eqref{Christoffel} that
\begin{multline}
\label{step1}
\lambda_n\big(ze^{\ic a};\mu\big) \leq |S_{\lfloor n/2m\rfloor}\big(ze^{\ic a},ze^{\ic a}\big)|^{-2}\int_\T |S_{\lfloor n/2m\rfloor}\big(ze^{\ic a},\eta\big)|^2\dd\mu(\eta) \\ \lesssim \frac1n\int_T |S_{\lfloor n/2m\rfloor}\big(ze^{\ic a},\eta\big)|^2\mu_n(\eta)|\dd\eta| + \frac1{n^2}\int_{\T\setminus T} |S_{\lfloor n/2m\rfloor}\big(ze^{\ic a},\eta\big)|^2\dd\mu(\eta),
\end{multline}
where the inequality on \( T \) follows from \cite[Equation~(4.25)]{And12}. It is known, see for example \cite[Lemma~2.1(ix)]{MastTot00}, that the doubling property is equivalent to
\begin{equation}
\label{doubling}
\mu_n(\eta) \lesssim (1+n|z-\eta|)^s\mu_n(z),
\end{equation}
uniformly \( z,\eta\in T \), where the parameter \( s \) depends only on the constant \( L \) in the doubling inequality \( \mu(2I) \leq L\mu(I) \). Choose \(m>(s+1)/2 \). Then \eqref{Sn-upper} and \eqref{doubling} yield that the first integral in \eqref{step1} can be estimates above by a constant times
\begin{multline}
\label{step2}
\frac{\mu_n(z)}n \int_T |S_{\lfloor n/2m\rfloor}\big(ze^{\ic a},\eta\big)|^2(1+n|z-\eta|)^s|\dd\eta| \\ \lesssim n\mu_n(z)\int_{T(z,\frac1n)}|\dd\eta| + \frac{\mu_n(z)}{n^{4m-s-1}}\int_{T\setminus T(z,\frac1n)}\frac{|\dd\eta|}{|z-\eta|^{4m-s}} \lesssim \mu_n(z).
\end{multline}
To estimate the second integral in \eqref{step1}, let us point out that \eqref{doubling} is a consequence of the inequality
\[
\int_I\dd\mu \lesssim \left(\frac{|I|+|J|+\dist(I,J)}{|J|}\right)^s\int_J\dd\mu, \quad I,J\subseteq T,
\]
where the constant is independent of \( I,J \), see \cite[Lemma~2.1(viii)]{MastTot00}. Therefore,
\begin{equation}
\label{decay}
\mu_n(z) \gtrsim n^{-s}, \quad z\in T.
\end{equation}
Thus, our choice of \( m \), \eqref{Sn-upper}, and \eqref{decay} imply that
\begin{equation}
\label{step3}
\frac1{n^2}\int_{\T\setminus T} |S_{\lfloor n/2m\rfloor}(z,\eta)|^2\dd\mu(\eta) \lesssim \frac1{n^{4m}}\int_{\T\setminus T}\frac{\dd\mu(\eta)}{|z-\eta|^{4m}} \lesssim \frac1{n^{4m}} \lesssim \frac1{n^{2s+2}} \lesssim \mu_n(z).
\end{equation}
The upper bound in \eqref{Christoffel-bounds} follows now by plugging estimates \eqref{step2} and \eqref{step3} into \eqref{step1}.

It only remains to prove the lower bound in \eqref{Christoffel-bounds}. Let \( z\in T^\prime \) and \( |a|\leq 2 \). Define
\[
Q_n(w) := w^{m(\lfloor n/2m\rfloor-1)}\int_\T S_{\lfloor n/2m\rfloor}(w,\eta)\big(n\mu_n(\eta)\big)^{1/2}|\dd\eta|,
\]
which is a polynomial of degree at most \( n-2m \). We get from \eqref{Sn-lower} and \eqref{doubling} that
\begin{eqnarray}
\big| Q_n\big(ze^{\ic a/n}\big) \big| & \gtrsim & \left|\int_{T(z,\frac1n)}\re\left(S_{\lfloor n/2m\rfloor}\big(ze^{\ic a/n},\eta\big)\right)(n\mu_n(\eta)\big)^{1/2}|\dd\eta|\right| \nonumber \\
& = & \left|\int_{-1}^1\frac1n\re\left(S_{\lfloor n/2m\rfloor}\big(e^{\ic a/n},e^{\ic t/n}\big)\right)(n\mu_n\big(ze^{\ic t/n}\big)\big)^{1/2}\dd t\right| \nonumber \\
\label{Qn-lower}
& \gtrsim &  \int_{-1}^1(n\mu_n\big(ze^{\ic t/n}\big)\big)^{1/2}\dd t \gtrsim (n\mu_n(z)\big)^{1/2}.
\end{eqnarray}
As \( S_n(w,\eta) \) is positive for \( w,\eta\in\T \), it follows from the normalization of \( S_n \) and Jensen's inequality that
\[
|Q_n(z)|^2 \leq \left(\int_T+\int_{\T\setminus T}\right)S_{\lfloor n/2m\rfloor}(z,\eta)n\mu_n(\eta)|\dd\eta|.
\]
Similarly to \eqref{step2}, the first integral above can be estimated as follows:
\[
\int_TS_{\lfloor n/2m\rfloor}(z,\eta)n\mu_n(\eta)|\dd\eta| \lesssim n\mu_n(z) \int_T(1+n|z-\eta|)^sS_{\lfloor n/2m\rfloor}(z,\eta)|\dd\eta|  \lesssim n\mu_n(z),
\]
where the first estimate follows from \eqref{doubling} and the second one from \eqref{Sn-upper}. Moreover, we also have that
\[
\int_{\T\setminus T}S_{\lfloor n/2m\rfloor}(z,\eta)n\mu_n(\eta)|\dd\eta| \lesssim n\int_{\T\setminus T}S_{\lfloor n/2m\rfloor}(z,\eta)|\dd\eta| \lesssim \frac1{n^{2(m-1)}} \leq \frac n{n^s} \lesssim n\mu_n(z),
\]
where we again used \eqref{Sn-upper} as well as \eqref{decay}. Altogether, we get that
\begin{equation}
\label{Qn-upper}
|Q_n(z)|^2 \lesssim n\mu_n(z), \quad z\in T^\prime.
\end{equation}

Let \( T_n \) be a polynomial of degree at most \( n-1 \) normalized to have value \( 1 \) at \( ze^{\ic a/n} \). It follows from \eqref{Qn-lower} and \eqref{Qn-upper} that
\begin{multline*}
\int|T_n|^2\dd\mu \geq \int_{T^\prime}|T_n|^2\dd\mu \gtrsim \int_{T^\prime}|T_n(\eta)|^2n\mu_n(\eta)|\dd\eta| \gtrsim \int_{T^\prime}|(T_nQ_n)(\eta)|^2|\dd\eta| \\ = \big|Q_n\big(ze^{\ic a/n}\big)\big|^2\int_{T^\prime}\frac{|(T_nQ_n)(\eta)|^2}{|Q_n(ze^{\ic a/n})|^2}|\dd\eta| \gtrsim n\mu_n(z)\int_{T^\prime}\frac{|(T_nQ_n)(\eta)|^2}{|Q_n(ze^{\ic a/n})|^2}|\dd\eta|,
\end{multline*}
where we also used \cite[Equation~(4.25)]{And12} for the second inequality. Since the polynomial in the last integral above is normalized to have value \( 1 \) at \( ze^{\ic a/n} \) and is of degree at most \( 2n \), \eqref{Christoffel} yields that
\begin{equation}
\label{Totik}
\lambda_n\big(ze^{\ic a/n};\mu\big) \gtrsim n\mu_n(z)\lambda_{2n}\big(ze^{\ic a/n};\sigma\big),
\end{equation}
where \( \sigma \) is the arclength measure on \( T^\prime \). Now, we know from \cite[Theorem~2.4]{Var13} that
\begin{equation}
\label{Varga}
n\lambda_{2n}(z;\sigma) = nK_{2n}^{-1}(z,z;\sigma) \gtrsim 1, \quad z\in T^\prime,
\end{equation}
with the constant independent of \( z \) and \( n \), where \( K_{2n}(z,w;\sigma) \) is defined exactly as in \eqref{kernels} with \( f_i(z) = \varphi_i(z;\sigma) \) and \( \varphi_i(z;\sigma) \) being \( i \)-th orthonormal polynomial with respect to \( \sigma \). Then since
\[
|K_{2n}(z,w;\sigma)| \leq \sqrt{K_{2n}(z,z;\sigma)}\sqrt{K_{2n}(w,w;\sigma)},
\]
which is simply an application of the Cauchy-Schwarz inequality, and by applying Bernstein-Walsh inequality to each variable of \( K_{2n}(z,w;\sigma) \) as it was done in \cite[Lemma~6.2]{LevLub07}, we get from \eqref{Varga}  that
\[
\big|K_{2n}\big(ue^{\ic b/n}, ve^{\ic c/n};\sigma\big)\big| \lesssim n, \quad u,v\in T^\prime, \quad |b|,|c|\leq 2.
\]
The lower bound in \eqref{Christoffel-bounds} now follows by restricting the above inequality to the diagonal and plugging it in \eqref{Totik}.

\subsection{Upper Estimate}

It follows from \eqref{expected-SV}, \eqref{OPUC-real-intensity}, and \eqref{symmetry} that
\[
\E(N_n) = 2\int_{-1}^1 \rho_n^{(1,0)}(x)\dd x \leq \frac2\pi\log n + \mathcal O(1) + 2\left(\int_{-1}^{-1+1/n}+\int_{1-1/n}^1\right)\rho_n^{(1,0)}(x)\dd x.
\]
We would like to show that the last two integrals are bounded above by an absolute constant.  We shall show this only for the integral on \( [1-\frac1n,1] \), the case of the other one being completely identical.  It holds that
\begin{equation}
\label{Ipm1}
\int_{1-1/n}^1 \rho_n^{(1,0)}(x)\dd x \leq \int_{1- 1/n}^1 \sqrt{\frac{K_n^{(1,1)}(x,x)}{K_n(x,x)}} \frac{\dd x}\pi = \int_0^1 \sqrt{\frac{K_n^{(1,1)}(1-\frac yn,1-\frac yn)}{n^2K_n(1-\frac yn,1-\frac yn)}} \frac{\dd y}\pi.
\end{equation}
which is an easy consequence of \eqref{real-intensity-Vanderbei}. As \( \mu \) is doubling in some neighborhood of \( 1 \), it follows from the Cauchy-Schwarz inequality, \eqref{Christoffel}, and the lower bound in \eqref{Christoffel-bounds} that
\[
\left|K_n\left(1+\frac un,1+\frac {\bar v}n\right)\right| \leq K_n^{1/2}\left(1+\frac un,1+\frac un\right) K_n^{1/2}\left(1+\frac {\bar v}n,1+\frac {\bar v}n\right) \lesssim \mu_n^{-1}(1)
\]
for \( |u|,|v|\leq 3/2 \). Consequently, the Cauchy integral formula for derivatives of holomorphic functions gives us
\[
\left|K_n^{(1,1)}\left(1+\frac un,1+\frac{\bar v}n\right)\right| = \left|\frac1{(2\pi\ic)^2}\int_{|\eta|=3/2}\int_{|\xi|=3/2}\frac{K_n(1+\frac \eta n,1+\frac {\bar \xi}n)}{(\frac\eta n-\frac un)^2(\frac{\bar\xi}n - \frac{\bar v}n)^2}\frac{\dd\eta}n\frac{\dd \bar\xi}n\right| \lesssim \frac{n^2}{\mu_n(1)}
\]
for  \( |u|,|v|\leq 1 \). The desired claim now follows from the above inequality combined with the upper estimate in \eqref{Christoffel-bounds}.

\subsection{Lower Estimate}

Under the current assumptions the measure \( \mu \) is doubling on the whole circle (see the explanation after the statement of the theorem) and therefore we only need to prove the lower estimate. To this end, observe that
\[
\E(N_n) > \frac2\pi\int_{-1+\frac{\log n}n}^{1-\frac{\log n}n}\frac{\sqrt{1-h_n^2(x)}}{1-x^2}\dd x \geq -\frac2\pi\sqrt{1-M_n^2} \log\left(\frac{\log n}n\right),
\]
where \( M_n \) is the maximum of \( |h_n(x)| \) on the interval of integration above. Thus, to prove \eqref{exact-asymp} it is enough to show that \( M_n = o(1) \) as \( n\to\infty \).

By the conditions of the theorem the sequence of the recurrence coefficients is absolutely summable. Hence, it follows from \cite[Theorem~1.5.3]{Simon1} that
\begin{equation}
\label{Normal}
1\lesssim |\Phi_n^*| \lesssim 1
\end{equation}
uniformly on \( \overline\D \). As \( |\Phi_n|=|\Phi_n^*| \) on \( \T \), it also follows from the Bernstein-Walsh inequality that
\begin{equation}
\label{BW}
|\Phi_n(z)| \lesssim |z|^n, \quad |z|\geq1.
\end{equation}
It can be readily checked using \eqref{recurrence} that
\[
\Phi_n^*(z) = 1 - z\sum_{k=0}^{n-1}\alpha_k\Phi_k(z).
\]
Therefore, we get from \eqref{BW} and the absolute summability of \( \alpha_k\)'s that
\begin{equation}
\label{ring}
|\Phi_n^*(z)| \lesssim \sum_{k=0}^{n-1}|\alpha_k||z|^{k+1} = \left(\sum_{k=0}^{m-1}+\sum_{k=m}^{n-1}\right)|\alpha_k||z|^{k+1} \lesssim |z|^m + \Lambda_m|z|^n
\end{equation}
for \( |z|\geq 1 \), where \( \Lambda_m := \sum_{k=m}^\infty|\alpha_k| \). That is, it holds that
\begin{equation}
\label{est-Phin}
|\Phi_n(z)| \lesssim \Lambda_m + |z|^{n-m}, \quad |z|\leq 1.
\end{equation}
Combining the above inequality with the lower bound in \eqref{Normal}, we see that
\begin{equation}
\label{est-bn}
|b_n(z)| \lesssim \Lambda_m + |z|^{n-m}, \quad |z|\leq 1.
\end{equation}
It further follows from the Cauchy integral formula for the derivatives that
\begin{equation}
\label{est-bn-prime1}
\big|\big(zb_{n-1}(z)\big)^\prime\big| \leq \int_{|\zeta|=r}\frac{|\zeta b_{n-1}(\zeta)|}{|\zeta-z|^2}\frac{|\mathrm d\zeta|}{2\pi} \leq r^2\frac{\Lambda_m+r^{n-m-1}}{r^2-|z|^2}, \quad |z|<r.
\end{equation}
On the other hand, the Bernstein inequality for polynomials on the disk of radius \( r\), \eqref{Normal}, and \eqref{est-Phin} yield that
\begin{equation}
\label{est-bn-prime2}
\max_{|z|\leq r}\big|\big(zb_{n-1}(z)\big)^\prime\big| \lesssim n\left(\Lambda_m + r^{n-m-1}\right).
\end{equation}
Now, take \( m \) to be the integer part of \( n/\log n \) and recall that \( \Lambda_m \lesssim m^{1-p} \) according to the condition placed on the recurrence coefficients. Thus, inequalities \eqref{est-bn-prime1} and \eqref{est-bn-prime2}, both applied with \( r = 1-\log n/n\), give
\begin{equation}
\label{est-sn}
\big|\big(zb_n(z)\big)^\prime\big| \lesssim \big(\log n\big)^{p-1} \left\{
\begin{array}{ll}
n^{3/2-p}, & |z| \leq 1-n^{-1/2}, \medskip \\
n^{2-p}, & 1-n^{-1/2}< |z| \leq 1-\frac{\log n}n.
\end{array}
\right.
\end{equation}
It also follows from \eqref{est-bn} that for \( x\in(-1,1) \) we have that
\[
\frac{1-x^2}{1-(xb_{n-1}(x))^2} \lesssim \left\{
\begin{array}{ll}
1, & |x| \leq 1-n^{-1/2}, \medskip \\
n^{-1/2}, & 1-n^{-1/2}< |x| \leq 1-\frac{\log n}n.
\end{array}
\right.
\]
Therefore, we deduce from \eqref{hn-rep} that \( M_n \lesssim (\log n)^{p-1}n^{3/2-p} = o(1) \) as desired.

\subsection{Proof of Proposition~\ref{prop:real-expected}}

Notice that the upper bound \eqref{upper-bound} remains valid in this case. We prove the lower bound as in the previous section by showing that the maximum of \( |h_{n+1}(x)| \) on \( [-1+\log n/n,1-\log n/n]\) behaves like \( o(1) \) as \( n\to\infty \). 

It was discovered by Geronimus \cite{Geronimus}, see also \cite{Wong09}, that
\begin{eqnarray*}
\Phi_n(z) &=& \Phi_n(z;\nu) - \frac{\Phi_n(1;\nu)K_n(z,1;\nu)}{t(1-t)^{-1}+K_n(1,1;\nu)} \\
&=& \Phi_n(z;\nu) - \beta_n\frac{\Phi_n^*(z;\nu)-\Phi_n(z;\nu)}{1-z},
\end{eqnarray*}
where we used Christoffel-Darboux formula \eqref{Kn} and set
\[
\beta_n := \frac{\varphi_n^2(1;\nu)}{t(1-t)^{-1}+K_n(1,1;\nu)}.
\]
Put \( s_n(z) := zb_n(z;\nu) \).
Then it holds that
\[
zb_n(z) = \frac{(1-z)s_n(z)-\beta_n(z-s_n(z))}{1-z(1+\beta_n)+\beta_ns_n(z)} = 1 - \frac{(1-z)(1-s_n(z))}{1-z(1+\beta_n)+\beta_ns_n(z)}.
\]
Therefore,
\[
1-\big(zb_n(z)\big)^2 = (1-z)^2(1-s_n(z))\frac{1+s_n(z)+2(s_n(z)-z)\frac{\beta_n}{1-z}}{(1-z(1+\beta_n)+\beta_ns_n(z))^2}
\]
and
\[
\big(zb_n(z)\big)^\prime = \frac{(1-z)^2(1+\beta_n)s_n^\prime(z)+\beta_n(2s_n(z)-1)}{(1-z(1+\beta_n)+\beta_ns_n(z))^2}.
\]
Thus, we get from \eqref{hn-rep} that
\begin{equation}
\label{G1}
h_{n+1}(z) = \frac{(1-z^2)(1+\beta_n)s_n^\prime(z)+\beta_n\frac{1+z}{1-z}(2s_n(z)-1)}{(1-s_n(z))(1+s_n(z)+\frac{2\beta_n}{1-z}(s_n(z)-z))}.
\end{equation}
It follows from the explanation given after the statement of the proposition that \( \beta_n \sim 1/n \) and therefore
\[
\beta_n\frac{1+x}{1-x}|2s_n(x)-1| \lesssim \frac{\beta_n}{1-x} \lesssim \frac1{\log n}, \quad -1 \leq x \leq 1-\frac{\log n}{n},
\]
where we used the fact that \( |s_n(z)|\leq 1 \) in \( \D \). It also follows from \eqref{est-sn} that
\[
(1-x^2)(1+\beta_n)|s_n^\prime(x)| \lesssim (\log n)^{p-1}n^{3/2-p}, \quad |x| \leq 1-\frac{\log n}n.
\]
Hence, the numerator of \eqref{G1} is of order \( o(1) \) as \( n\to\infty \) on the interval of interest. Similarly, we see from \eqref{est-bn} that the denominator behaves like \( 1 + o(1) \) there, which finishes the proof of the proposition.

\section{Proof of Theorem~\ref{thm:complex-intensity}}

Let us modify expression \eqref{complex-intensity-Vanderbei} for \( \rho_n^{(0,1)} \) to make it more amenable to the asymptotic analysis. Write
\begin{equation}
\label{repr1}
\pi\left( K_n(z,z)^2-|K_n(z,\overline z)|^2\right)^{3/2}\rho_n^{(0,1)}(z) := S_1(z) + S_2(z) + S_3(z),
\end{equation}
where the functions \( S_i(z) \) correspond to the respective summands in \eqref{complex-intensity-Vanderbei}. For brevity, put
\[
S_n(z,w) := (\varphi_n^*)^\prime(z)\varphi_n^*(\overline w)-\varphi_n^\prime(z)\varphi_n(\overline w).
\]
Then it follows from \eqref{Kn}--\eqref{Kn11} that \( S_1(z) \) is equal to
\begin{eqnarray}
\label{S11}
&  & \frac{1+|z|^2}{(1-|z|^2)^2}K_n(z,z)\left(K_n(z,z)^2 - |K_n(z,\overline z)|^2\right) \\
\label{S12}
& + & 2\frac{\re(zS_n(z,z))}{(1-|z|^2)^2}\left(K_n(z,z)^2 - |K_n(z,\overline z)|^2\right) \\
\label{S13}
& + & \frac{|(\varphi_n^*)^\prime(z)|^2-|\varphi_n^\prime(z)|^2}{1-|z|^2}\left(K_n(z,z)^2 - |K_n(z,\overline z)|^2\right),
\end{eqnarray}
\( S_2(z) \) is equal to
\begin{eqnarray}
\label{S21}
&  & -|z|^2K_n(z,z)\left(\frac{K_n(z,z)^2}{(1-|z|^2)^2} - \frac{|K_n(z,\overline z)|^2}{|1-z^2|^2}\right) \\
\label{S22}
& - & 2\frac{K_n(z,z)^2\re(zS_n(z,z))}{(1-|z|^2)^2} - 2\frac{K_n(z,z)\re\big(zK_n(z,\overline z)\overline{S_n(z,\overline z)}\big)}{|1-z^2|^2} \\
\label{S23}
& - & \frac{K_n(z,z)|S_n(z,z)|^2}{(1-|z|^2)^2} - \frac{K_n(z,z)|S_n(z,\overline z)|^2}{|1-z^2|^2},
\end{eqnarray}
and \( S_3(z) \) is equal to
\begin{eqnarray}
\label{S31}
&  &  \frac{K_n(z,z)|K_n(z,\overline z)|^2}{1-|z|^2}\left(\frac{1-|z|^4}{|1-z^2|^2}-1\right) \\
\label{S32}
& + & 2\frac{|K_n(z,\overline z)|^2}{1-|z|^2}\re\left(\frac{\overline zS_n(z,z)}{1-\overline z^2}\right) + 2\frac{K_n(z,z)}{1-|z|^2}\re\left(\frac{\overline zK_n(z,\overline z)\overline{S_n(z,\overline z)}}{1-\overline z^2}\right) \\
\label{S33}
& + & \frac2{1-|z|^2}\re\left(\frac{K_n(z,\overline z)S_n(z,z)\overline{S_n(z,\overline z)}}{1-\overline z^2}\right),
\end{eqnarray}
where we used the indentity \( 2\re(z^2) = 1 + |z|^4-|1-z^2|^2 \) in \eqref{S31}. Then we can rewrite \eqref{repr1} as
\begin{equation}
\label{repr2}
\pi\left( K_n(z,z)^2-|K_n(z,\overline z)|^2\right)^{3/2}\rho_n^{(0,1)}(z) := \Sigma_{n,1}(z) - \Sigma_{n,2}(z) + \Sigma_{n,3}(z),
\end{equation}
where \( \Sigma_{n,1} \) is the sum of \eqref{S11}, \eqref{S21}, and \eqref{S31}, \( \Sigma_{n,2} \) is the sum of \eqref{S12}, \eqref{S22}, and \eqref{S32}, and \( \Sigma_{n,3} \) is the sum of \eqref{S13}, \eqref{S23}, and \eqref{S33}. One can readily verify that
\begin{equation}
\label{Sigma1}
\Sigma_{n,1}(z) = \frac{K_n(z,z)^3}{(1-|z|^2)^2} - K_n(z,z)|K_n(z,\overline z)|^2\left(\frac2{(1-|z|^2)^2} - \frac1{|1-z^2|^2}\right).
\end{equation}
Furthermore, one can check that the sum of \eqref{S12} and the first summands of \eqref{S22} and  \eqref{S32} is equal to
\[
2 \frac{\re\big((\overline z-z)\overline{K_n(z,\overline z)}(\varphi_n^*(z)^2-\varphi_n(z)^2)S_n(z,z)\big)}{(1-|z|^2)^2|1-z^2|^2},
\]
while the sum of the last summands of \eqref{S22} and \eqref{S32} is equal to
\[
2 \frac{\re\big((z-\overline z)\overline{K_n(z,\overline z)}(|\varphi_n^*(z)|^2-|\varphi_n(z)|^2)S_n(z,\overline z)\big)}{(1-|z|^2)^2|1-z^2|^2}.
\]
By adding up the last two expressions and simplifying, we get that
\begin{equation}
\label{Sigma2}
\Sigma_{n,2}(z) = \frac{8\im(z)\im\big(\varphi_n^*(z)\overline{\varphi_n(z)}\big)\re((\varphi_n^*)^\prime(z)\varphi_n(z)-\varphi_n^\prime(z)\varphi_n^*(z))}{(1-|z|^2)^2|1-z^2|^2}.
\end{equation}
To compute \( \Sigma_{n,3} \), notice that the sum of the first summands of \eqref{S13} and \eqref{S23} is equal to
\[
-\frac{K_n(z,z)|(\varphi_n^*)^{\prime}(z)\varphi_n(z)-\varphi_n^{\prime}(z)\varphi_n^*(z)|^2}{(1-|z|^2)^2}.
\]
The remaining summand of \eqref{S13} is equal to
\[
\big(|\varphi_n^\prime(z)|^2-  |(\varphi_n^*)^\prime(z)|^2\big) \frac{|\varphi_n(z)|^4+|\varphi_n^*(z)|^4 - 2\re\big(\varphi_n^*(z)^2\overline{\varphi_n(z)^2}\big)}{(1-|z|^2)|1-z^2|^2},
\]
while the remaining summand of \eqref{S23} is equal to
\begin{multline*}
|\varphi_n^\prime(z)|^2\frac{|\varphi_n(z)|^4-|\varphi_n(z)\varphi_n^*(z)|^2}{(1-|z|^2)|1-z^2|^2} - |(\varphi_n^*)^\prime(z)|^2\frac{|\varphi_n^*(z)|^4-|\varphi_n(z)\varphi_n^*(z)|^2}{(1-|z|^2)|1-z^2|^2} \\ + 2K_n(z,z)\frac{\re\big((\varphi_n^*)^\prime(z)\varphi_n^*(z)\overline{\varphi_n^\prime(z)\varphi_n(z)\big)}}{|1-z^2|^2}.
\end{multline*}
Moreover, \eqref{S33} can be rewritten as
\begin{multline*}
|\varphi_n^\prime(z)|^2\frac{2\re\big(\varphi_n^*(z)^2\overline{\varphi_n(z)^2}\big)-2|\varphi_n(z)|^4}{(1-|z|^2)|1-z^2|^2} - |(\varphi_n^*)^\prime(z)|^2\frac{2\re\big(\varphi_n^*(z)^2\overline{\varphi_n(z)^2}\big)-2|\varphi_n^*(z)|^4}{(1-|z|^2)|1-z^2|^2} \\ - 2K_n(z,z)\frac{\re\big((\varphi_n^*)^\prime(z)\varphi_n^*(z)\overline{\varphi_n^\prime(z)\varphi_n(z)}+(\varphi_n^*)^\prime(z)\varphi_n(z)\overline{\varphi_n^\prime(z)\varphi_n^*(z)}\big)}{|1-z^2|^2}.
\end{multline*}
By adding the last four expressions together we get that
\begin{equation}
\label{Sigma3}
\Sigma_{n,3}(z) = K_n(z,z)|(\varphi_n^*)^\prime(z)\varphi_n(z)-\varphi_n^\prime(z)\varphi_n^*(z)|^2\left(\frac1{|1-z^2|^2} - \frac1{(1-|z|^2)^2} \right).
\end{equation}
Notice that
\begin{equation}
\label{derivative}
\frac{(\varphi_n^*)^\prime(z)\varphi_n(z)-\varphi_n^\prime(z)\varphi_n^*(z)}{\phi_n^2(z)} = \left\{
\begin{array}{ll}
-b_n^\prime(z), & |z|<1, \medskip \\
(b_n^{-1})^\prime(z), & |z|>1,
\end{array}
\right.
\end{equation}
where
\[
\phi_n(z) := \left\{
\begin{array}{ll}
\varphi_n^*(z), & |z|<1, \medskip \\
\varphi_n(z), & |z|>1.
\end{array}
\right.
\]
Finally, the assumption \( \alpha_i \to 0 \) as \( i\to\infty \) implies that
\begin{equation}
\label{limitto0}
\left\{
\begin{array}{ll}
b_n(z) \to 0, & \text{locally uniformly in } |z|<1, \medskip \\
b_n^{-1}(z)\to 0, & \text{locally uniformly in } |z|>1,
\end{array}
\right.
\end{equation}
as \( n\to\infty \) according to \cite[Theorem 1.7.4]{Simon1} and since \( b_n^{-1}(z)=b_n(1/z) \). By recalling \eqref{Kn} and pugging  \eqref{derivative} into \eqref{Sigma2}, \eqref{Sigma3} and using \eqref{limitto0}, we get that
\begin{equation}
\label{almostthere1}
\frac{(\Sigma_{n,2}+\Sigma_{n,3})(z)}{|\phi_n(z)|^6} \to 0
\end{equation}
as \( n\to\infty \) locally uniformly in \( \overline\C\setminus\T \). Similarly, we get that
\begin{equation}
\label{almostthere2}
\frac{\Sigma_{n,1}(z)}{|\phi_n(z)|^6} \to \frac1{|1-|z|^2|}\left(\frac1{(1-|z|^2)^2} - \frac1{|1-z^2|^2}\right)^2
\end{equation}
as \( n\to\infty \) locally uniformly in \( \overline\C\setminus\T \). Finally, since
\[
\frac{K_n(z,z)^2 - |K_n(z,\overline z)|^2}{|\phi_n(z)|^4} \to \left(\frac1{(1-|z|^2)^2} - \frac1{|1-z^2|^2}\right)
\]
as \( n\to\infty \) locally uniformly in \( \overline\C\setminus\T \), the claim of the theorem follows from \eqref{almostthere2}, \eqref{almostthere1}, and \eqref{repr2}.

\section{Proof of Theorem~\ref{thm:complex-expected}}

Given an arc \( S \), it follows from \cite[Theorem 3.1]{PritYgr15} that
\begin{multline}
\label{PrY}
\E \left[ \left| \frac1nN_n(\Omega(S,\delta))-\frac{|S|}{2\pi}  \right| \right] \lesssim \\
\frac1{\sqrt n\delta}\left[\log \sum_{k=0}^{n-1}\E \big[|\eta_k| \big] + \max_{0\leq k \leq n-1}\log \|\varphi_k\|_{\infty} -\frac12\E \big[ \log |D_{n-1}| \big]  \right]^{1/2}
\end{multline}
with some absolute constant, where \( D_n: =\eta_n\kappa_n\sum_{k=0}^n\eta_k a_{0,k} \). Since the random variables \( \{\eta_k\} \) are standard Gaussian, we have that the first expression in square brackets is of the order \( \log n \). Furthermore, it follows from \cite[Equation~(1.5.17)]{Simon1}, the Cauchy-Schwarz inequality, and \eqref{leading-coeff} that
\[
\log\|\varphi_k\|_\infty \leq \log\kappa_k + \sum_{j=0}^{k-1}|\alpha_j| \leq  \log\kappa_n + \sum_{j=0}^{n-1}|\alpha_j| \leq \log\kappa_n + \left(n\sum_{j=0}^{n-1}\alpha_j^2\right)^{1/2}.
\]
Thus, the middle term on the right-hand side of \eqref{PrY} can be estimated above by
\[
\log\kappa_n + \left(n\sum_{j=0}^{n-1}\log\frac1{1-\alpha_j^2}\right)^{1/2} = \log\kappa_n  + (2n\log\kappa_n)^{1/2} \lesssim n\sqrt{\varepsilon_n},
\]
where we used \eqref{leading-coeff} and \eqref{Reg}. Finally, since there is constant \( C \) such that \( \E\big[\log |\eta_0+z|\big]\geq C \) for all \( z\in\C \), it holds that
\[
\E \big[ \log |D_n| \big] = \log \kappa_n + \E\big[\log |\eta_n|\big] + \E\left[\log \left|\eta_0 + \sum_{i=1}^{n-1}\eta_ia_{0,i}\right|\right] \geq 2C,
\]
(recall also that \(\kappa_n\geq1 \) and \( a_{0,0} = 1 \) as \( \mu \) is a probability measure). Plugging the above estimates into \eqref{PrY} proves the theorem.

\section{Proof of Theorem~\ref{thm:complex-expected2}}

It follows from \eqref{expected-SV} that
\begin{eqnarray*}
\frac1n \E\big[N_n\big(\Omega(S,\tau_1,\tau_2)\big)\big] & = & \frac1n\iint_{\Omega(S,\tau_1,\tau_2)}\rho_n^{(0,1)}(z)\dd A \\
& = & \frac1{2n^2}\int_S\int_{\tau_1}^{\tau_2}\rho_n^{(0,1)}\left(z\left(1+\frac\tau{2n}\right)\right)\left(1+\frac\tau{2n}\right)\dd\tau|\dd z|.
\end{eqnarray*}
Hence, it is enough to show that
\begin{equation}
\label{needed}
\frac1{2n^2}\rho_n^{(0,1)}\left(z\left(1+\frac\tau{2n}\right)\right) \to \frac1{2\pi}\left(\frac{H^\prime(\tau)}{H(\tau)}\right)^\prime
\end{equation}
uniformly for \( z\in S \) and \( \tau \) on compact subsets of the real line.

Since \( \alpha_k\to 0 \) as \( k \to\infty \), the measure \( \mu \) is regular in the sense of Ullman-Stahl-Totik, see \eqref{leading-coeff} and \eqref{Reg}. Therefore, \cite[Theorem~6.3]{LevLub07} is applicable on \( S \). That is, it holds that
\begin{equation}
\label{scaling-limit-1}
\lim_{n\to\infty} K_n(z_{n,u},z_{n,\overline v})K_n^{-1}(z,z) = H(u+v)
\end{equation}
uniformly for \( z \in S \) and \( u,v \) on compact subsets of \( \C \), where \( z_{n,a}:=z(1+a/n) \). Moreover, we have that
\[
\frac{\partial^{i+j}}{\partial u^i\partial v^j} K_n(z_{n,u},z_{n,\overline v}) = \frac{z^{i-j}}{n^{i+j}} K_n^{(i,j)}(z_{n,u},z_{n,\overline v})
\]
for any non-negative integers \( i,j \). Thus, it follows from Cauchy's integral formula and \eqref{scaling-limit-1} that
\begin{equation}
\label{scaling-limit-2}
\lim_{n\to\infty} \frac{z^{i-j}}{n^{i+j}}\frac{K_n^{(i,j)}(z_{n,u},z_{n,\overline v})}{K_n(z,z)} = H^{(i+j)}(u+v)
\end{equation}
uniformly for \( z \in S \) and \( u,v \) on compact subsets of \( \C \).

In another connection, since \( \alpha_i\to0 \) as \( i\to\infty \), \cite[Theorem~4]{MaNevTot87} states that
\[
\lim_{n\to\infty} \max_{z\in\T}|\varphi_n(z)|^2K_n^{-1}(z,z) = 0.
\]
By compactness, the set \( S \) can be covered by finitely many closed subarcs \( I_j\subset\T\setminus\{\pm1\} \) such that \( \mu^\prime \) is continuous and positive on each \( I_j \). Since \( \cup_j I_j \) is separated from \( \pm1 \), the Christoffel-Darboux formula \eqref{Kn} and the above limit yield that
\[
\lim_{n\to\infty} \max_{z\in \cup_j I_j}K_n(z,\overline z)K_n^{-1}(z,z) = 0.
\]
Since \( K_n(z,\overline z) \) is a polynomial of degree \( 2n-2 \), it follows from the Bernstein-Walsh inequality (see for example \cite[Lemma~6.1]{LevLub07}), applied on each \( I_j \) separately, that
\[
|K_n(z_{n,a},\overline z_{n,a})| \lesssim \max_{z\in \cup_jI_j}|K_n(z,\overline z)|
\]
uniformly in \( n \) and \( a \) on compact subsets of \( \C \). Thus, it holds that
\begin{equation}
\label{scaling-limit-3}
\lim_{n\to\infty} |K_n(z_{n,a},\overline z_{n,a})|K_n^{-1}(z,z) =0
\end{equation}
uniformly for \( z \in S \) and \( a \) on compact subsets of \( \C \). Moreover, since
\[
K_n^{(1,0)}(z_{n,a},\overline z_{n,a}) = \frac n{2z}\frac\partial{\partial a} K_n(z_{n,a},\overline z_{n,a}),
\]
it follows from \eqref{scaling-limit-3} and Cauchy's integral formula that
\begin{equation}
\label{scaling-limit-4}
\lim_{n\to\infty} n^{-1}|K_n^{(1,0)}(z_{n,a},\overline z_{n,a})|K_n^{-1}(z,z) =0
\end{equation}
uniformly for \( z \in S \) and \( a \) on compact subsets of \( \C \).

The desired claim \eqref{needed} now is an immediate consequence of \eqref{complex-intensity-Vanderbei} and \eqref{scaling-limit-2}--\eqref{scaling-limit-4}.

\section{Acknowledgments}

The second author would like to thank his Ph.D. advisor Igor Pritsker for his many helpful discussions concerning this project.

\end{document}